# LOWER LIMITS AND EQUIVALENCES
# FOR CONVOLUTION TAILS

By Serguei Foss[1,2] and Dmitry Korshunov[2]

*Heriot–Watt University and Sobolev Institute of Mathematics*

Suppose $F$ is a distribution on the half-line $[0,\infty)$. We study the limits of the ratios of tails $\overline{F * F}(x)/\overline{F}(x)$ as $x \to \infty$. We also discuss the classes of distributions $\mathcal{S}$, $\mathcal{S}(\gamma)$ and $\mathcal{S}^*$.

**1. Introduction.** Let $F$ be a distribution on the half-line $[0,\infty)$ with unbounded support, that is, $\overline{F}(x) \equiv F(x,\infty) > 0$ for any $x$. Let $a \in (0,\infty]$ be the mean value of $F$. By the Laplace transform of $F$ at the point $\gamma \in \mathbf{R}$ we mean

$$\varphi(\gamma) = \int_0^\infty e^{\gamma x} F(dx) \in (0,\infty].$$

Put

$$\hat{\gamma} = \sup\{\gamma : \varphi(\gamma) < \infty\} \in [0,\infty].$$

Note that the function $\varphi(\gamma)$ is monotone continuous in the interval $(-\infty, \hat{\gamma})$, and $\varphi(\hat{\gamma}) = \lim_{\gamma \uparrow \hat{\gamma}} \varphi(\gamma) \in [1,\infty]$.

We distinguish all the distributions on $[0,\infty)$ according to the value of $\hat{\gamma}$. If $\hat{\gamma} = 0$, then we say that the distribution $F$ is *heavy-tailed*; in that case $\varphi(\gamma) = \infty$ for any $\gamma > 0$. If $\hat{\gamma} > 0$, then we call the distribution $F$ *light-tailed*; this happens if and only if, for some $\gamma > 0$, $\overline{F}(x) = o(e^{-\gamma x})$ as $x \to \infty$.

The main results of this paper are the following Theorems 1, 2 and 3 which relate the tail behavior of the convolution $F * F$ to that of $F$.

Theorem 1. *For any heavy-tailed distribution $F$,*

$$\liminf_{x \to \infty} \frac{\overline{F * F}(x)}{\overline{F}(x)} = 2.$$

Received July 2005; revised March 2006.

[1]Supported by EPSRC Travel Grant GR/T27099/01.

[2]Supported by Science Foundation Ireland Grant SFI 04/RP1/I512.

We were motivated by the nice paper of Rudin [12]. Theorem $2^*$ of that paper states that, for any independent stopping time $\tau$, the distribution tail of the sum $S_\tau = \xi_1 + \cdots + \xi_\tau$ of i.i.d. random variables with common distribution $F$ satisfies the relation

$$(1) \qquad \liminf_{x \to \infty} \frac{\mathbf{P}\{S_\tau > x\}}{\overline{F}(x)} = \mathbf{E}\tau,$$

provided (i) $\mathbf{E}\xi^p = \infty$ and (ii) $\mathbf{E}\tau^p < \infty$ for some positive integer $p$. Unfortunately, condition (i) rules out a lot of distributions of interest, say, in the theory of subexponential distributions. For example, log-normal and Weibull-type distributions do not satisfy (i). Theorem 1 is restricted to the case $\tau = 2$, but here extends Rudin's result to the class of all heavy-tailed distributions. The reasons for the restriction to $\tau = 2$ come from the proof of Theorem 1 but in fact are rather deep; we provide more detailed comments in Section 4 which is devoted to the proof. Note that the case $\tau = 2$ is of genuine interest in itself.

The counterpart of Theorem 1 in the light-tailed case is stated next.

THEOREM 2. *Let $\hat\gamma \in (0, \infty]$, so that $\varphi(\hat\gamma) \in (1, \infty]$. If, for any fixed $y > 0$,*

$$(2) \qquad \liminf_{x \to \infty} \frac{\overline{F}(x-y)}{\overline{F}(x)} \geq e^{\hat\gamma y},$$

*then*

$$\liminf_{x \to \infty} \frac{\overline{F * F}(x)}{\overline{F}(x)} = 2\varphi(\hat\gamma).$$

The proof follows from Lemmas 7, 8 and 9 in Section 6. It turns out that condition (2) is essential for the conclusion of Theorem 2 to hold; see the counterexamples in Section 9. Note also that condition (2) weakens the commonly used assumption that $F \in \mathcal{L}(\hat\gamma)$; see Section 8. We give here a few simple examples where conditions of Theorem 2 *are* satisfied. If $\overline{F}(x) \sim c_1 e^{-c_2 x^2}$, then $\hat\gamma = \infty$, $\varphi(\hat\gamma) = \infty$, the condition (2) is met, and we have $\overline{F * F}(x)/\overline{F}(x) \to \infty$ as $x \to \infty$. If $F$ is the exponential distribution with parameter $\alpha$, then $\hat\gamma = \alpha$, $\varphi(\hat\gamma) = \infty$, the condition (2) is met, and we again have the convergence $\overline{F * F}(x)/\overline{F}(x) \to \infty$. If $\overline{F}(x) = l(x)e^{-\alpha x}$ where $l(x)$ is positive, slowly varying at infinity, and integrable, then $\hat\gamma = \alpha$ and $\varphi(\hat\gamma) < \infty$.

The third theorem may be considered as the final point in a chain of results in this direction; see [2, 3, 4, 7, 10, 11, 13] and the references therein.

THEOREM 3. *Let $F$ be any distribution on $[0, \infty)$ with unbounded support. Assume that*

$$\frac{\overline{F * F}(x)}{\overline{F}(x)} \to c \qquad as\ x \to \infty,$$

*where $c \in (0, \infty]$. Then $c = 2\varphi(\hat{\gamma})$.*

The proof of Theorem 3 is given in Section 7. In our proof we use real analytic and direct probabilistic methods; no Banach algebra technique is involved. We would like to emphasize that we do not assume that $F \in \mathcal{L}(\hat{\gamma})$ or that $c$ is finite.

In Sections 3 and 5 we study local properties of convolutions. Section 8 discusses the classes of distributions $\mathcal{S}$, $\mathcal{S}(\gamma)$ and $\mathcal{S}^*$ from the point of view of the results obtained. In Section 10 we consider briefly the convolution of nonidentical distributions.

**2. Characterization of heavy-tailed distributions.** In the sequel we need the following existence result which generalizes a lemma by Rudin ([12], page 989) onto the whole class of heavy-tailed distributions. Fix any $\delta \in (0, 1]$.

LEMMA 1. *If a random variable $\xi \geq 0$ has a heavy-tailed distribution, then there exists a function $h : \mathbf{R}^+ \to \mathbf{R}^+$ such that:*

(i) *$h$ is subadditive, that is, $h(x) \leq h(y) + h(x - y)$ for any $0 \leq y \leq x$;*
(ii) *$h(x) = o(x)$ as $x \to \infty$;*
(iii) *$\mathbf{E}e^{h(\xi)} \leq 1 + \delta$;*
(iv) *$\mathbf{E}\xi e^{h(\xi)} = \infty$.*

PROOF. The proof is a straightforward modification of the corresponding lemma of [12]. Put $x_0 = 0$. Choose $x_1 \geq 2$ so that $\overline{F}(x_1) < \delta/2$. Choose $\varepsilon_1 > 0$ so that

$$\mathbf{E}\{e^{\varepsilon_1 \xi}; \xi \leq x_1\} < 1.$$

By induction we construct an increasing sequence $x_n$ and a decreasing sequence $\varepsilon_n > 0$ such that $x_n \geq 2^n$ and $\overline{F}(x_n) < \delta/2^n$ for any $n \geq 1$, and

$$\mathbf{E}\{e^{\varepsilon_n \xi}; \xi \in (x_{n-1}, x_n]\} = \delta/2^{n-1} \qquad \text{for any } n \geq 2.$$

For $n = 1$ this is already done. Make the induction hypothesis for some $n \geq 2$. Due to heavy-tailedness, there exists $x_{n+1} \geq 2^{n+1}$ so large that $\overline{F}(x_{n+1}) < \delta/2^{n+1}$ and

$$\mathbf{E}\{e^{\varepsilon_n \xi}; \xi \in (x_n, x_{n+1}]\} \geq \delta.$$

Since $\overline{F}(x_n) < \delta/2^n$, we can find $\varepsilon_{n+1} \in (0, \varepsilon_n)$ such that

$$\mathbf{E}\{e^{\varepsilon_{n+1} \xi}; \xi \in (x_n, x_{n+1}]\} = \delta/2^n.$$

Our induction hypothesis now holds with $n + 1$ in place of $n$ as required.

Define $h(0) = 0$, $h(x) = \varepsilon_n x$ if $x \in (x_{n-1}, x_n]$, $n \geq 1$. This function is subadditive. Indeed, for any $0 \leq y \leq x$, we have $y \in (x_i, x_{i+1}]$, $x - y \in (x_j, x_{j+1}]$ and $x \in (x_n, x_{n+1}]$ for some $i$, $j$ and $n$ where $i$, $j \leq n$. Then

$$h(y) + h(x - y) = \varepsilon_i y + \varepsilon_j (x - y) \geq \varepsilon_n y + \varepsilon_n (x - y) = \varepsilon_n x = h(x),$$

due to the monotonicity of $\varepsilon_n$. Note that the function $h(x)$ is not monotone.

Next,

$$\mathbf{E} e^{h(\xi)} = \mathbf{E}\{e^{\varepsilon_1 \xi}; \xi \leq x_1\} + \sum_{n=2}^{\infty} \mathbf{E}\{e^{\varepsilon_n \xi}; \xi \in (x_{n-1}, x_n]\}$$

$$\leq 1 + \sum_{n=2}^{\infty} \delta/2^{n-1} = 1 + \delta.$$

On the other hand,

$$\mathbf{E} \xi e^{h(\xi)} = \sum_{n=1}^{\infty} \mathbf{E}\{\xi e^{\varepsilon_n \xi}; \xi \in (x_{n-1}, x_n]\}$$

$$\geq \sum_{n=2}^{\infty} x_{n-1} \mathbf{E}\{e^{\varepsilon_n \xi}; \xi \in (x_{n-1}, x_n]\}$$

$$\geq \sum_{n=2}^{\infty} 2^{n-1} \delta/2^{n-1} = \infty.$$

Note also that necessarily $\lim_{n \to \infty} \varepsilon_n = 0$ [otherwise $\liminf_{x \to \infty} h(x)/x > 0$ and $\xi$ is light-tailed], and (ii) follows. The proof of the lemma is complete. $\square$

## 3. Heavy tails: local properties.
Lemma 1 provides a useful tool for proving upper bounds in lower limit assertions for convolution of densities.

Let $\mu$ be a $\sigma$-finite measure on $[0, \infty)$. Suppose that the distribution $F$ on $[0, \infty)$ is absolutely continuous with respect to $\mu$, and let $f(x)$ be the corresponding density, that is, the Radon–Nikodym derivative of $F$ with respect to $\mu$. Consider i.i.d. random variables $\xi_1$, $\xi_2, \ldots$ with common distribution $F$ and an independent stopping time $\tau$. Assume that the distribution of the randomly stopped sum $S_\tau = \xi_1 + \cdots + \xi_\tau$ has a density $f^{*\tau}(x)$ with respect to $\mu$. In the present section we are interested in the limiting behavior of the ratio $\frac{f^{*\tau}(x)}{f(x)}$ as $x \to \infty$; more precisely, the next lemma deals with the "$\liminf$" for that ratio and generalizes Theorem 4 of [12].

LEMMA 2. *If $F$ is heavy-tailed, then*

$$\liminf_{x \to \infty} \frac{f^{*\tau}(x)}{f(x)} \leq \mathbf{E} \tau,$$

*provided $\tau$ is light-tailed, that is, $\mathbf{E} e^{\kappa \tau} < \infty$ for some $\kappa > 0$.*

Lemma 2 implies a corollary for "local subexponentiality" which can be found in Section 8.

PROOF OF LEMMA 2. Assume the contrary, that is, there exist $\delta > 0$ and $x_0$ such that

$$(3) \qquad f^{*\tau}(x) \geq (\mathbf{E}\tau + \delta)f(x) \qquad \text{for all } x > x_0.$$

Since $\mathbf{E}e^{\kappa\tau} < \infty$,

$$(4) \qquad \mathbf{E}\tau(1 + \varepsilon)^{\tau-1} \leq \mathbf{E}\tau + \delta/2$$

for some sufficiently small $\varepsilon > 0$. By Lemma 1, there exists a subadditive function $h(x) > 0$ such that $\mathbf{E}e^{h(\xi_1)} \leq 1 + \varepsilon$ and $\mathbf{E}\xi_1 e^{h(\xi_1)} = \infty$. For any random variable $\zeta$ and positive $t$, put $\zeta^{[t]} = \min\{\zeta, t\}$. Then

$$\frac{\mathbf{E}(\xi_1^{[t]} + \cdots + \xi_\tau^{[t]})e^{h(\xi_1 + \cdots + \xi_\tau)}}{\mathbf{E}\xi_1^{[t]}e^{h(\xi_1)}}$$

$$= \sum_{n=1}^{\infty} \frac{\mathbf{E}(\xi_1^{[t]} + \cdots + \xi_n^{[t]})e^{h(\xi_1 + \cdots + \xi_n)}}{\mathbf{E}\xi_1^{[t]}e^{h(\xi_1)}}\mathbf{P}\{\tau = n\}$$

$$= \sum_{n=1}^{\infty} n\frac{\mathbf{E}\xi_1^{[t]}e^{h(\xi_1 + \cdots + \xi_n)}}{\mathbf{E}\xi_1^{[t]}e^{h(\xi_1)}}\mathbf{P}\{\tau = n\}$$

$$\leq \sum_{n=1}^{\infty} n\frac{\mathbf{E}\xi_1^{[t]}e^{h(\xi_1) + \cdots + h(\xi_n)}}{\mathbf{E}\xi_1^{[t]}e^{h(\xi_1)}}\mathbf{P}\{\tau = n\},$$

by subadditivity. Hence,

$$\frac{\mathbf{E}(\xi_1^{[t]} + \cdots + \xi_\tau^{[t]})e^{h(\xi_1 + \cdots + \xi_\tau)}}{\mathbf{E}\xi_1^{[t]}e^{h(\xi_1)}}$$

$$(5) \qquad \begin{aligned} &\leq \sum_{n=1}^{\infty} n\frac{\mathbf{E}\xi_1^{[t]}e^{h(\xi_1)}(\mathbf{E}e^{h(\xi_2)})^{n-1}}{\mathbf{E}\xi_1^{[t]}e^{h(\xi_1)}}\mathbf{P}\{\tau = n\} \\ &\leq \sum_{n=1}^{\infty} n(1 + \varepsilon)^{n-1}\mathbf{P}\{\tau = n\} \\ &\leq \mathbf{E}\tau + \delta/2, \end{aligned}$$

by (4). On the other hand, since $(\xi_1 + \cdots + \xi_\tau)^{[t]} \leq \xi_1^{[t]} + \cdots + \xi_\tau^{[t]}$,

$$\frac{\mathbf{E}(\xi_1^{[t]} + \cdots + \xi_\tau^{[t]})e^{h(\xi_1 + \cdots + \xi_\tau)}}{\mathbf{E}\xi_1^{[t]}e^{h(\xi_1)}}$$

$$(6) \qquad \geq \frac{\mathbf{E}(\xi_1 + \cdots + \xi_\tau)^{[t]} e^{h(\xi_1 + \cdots + \xi_\tau)}}{\mathbf{E}\xi_1^{[t]} e^{h(\xi_1)}}$$

$$= \frac{\int_0^\infty x^{[t]} e^{h(x)} f^{*\tau}(x)\mu(dx)}{\int_0^\infty x^{[t]} e^{h(x)} f(x)\mu(dx)}.$$

Since $\mathbf{E}\xi_1 e^{h(\xi_1)} = \infty$,

$$\int_0^\infty x^{[t]} e^{h(x)} f(x)\mu(dx) \to \infty \qquad \text{as } t \to \infty.$$

Therefore, it follows from (3) that

$$\liminf_{t\to\infty} \frac{\int_0^\infty x^{[t]} e^{h(x)} f^{*\tau}(x)\mu(dx)}{\int_0^\infty x^{[t]} e^{h(x)} f(x)\mu(dx)} \geq \mathbf{E}\tau + \delta.$$

Substituting this into (6), we get a contradiction to (5) for sufficiently large $t$. The proof is complete. $\square$

**4. Heavy tails: proof of Theorem 1.** First we restate in Lemma 3 below Theorem 1* of [12] in terms of probability distributions and stopping times. This result also follows immediately from the inequality (20).

LEMMA 3. *For any distribution $F$ on $[0, \infty)$ and any independent stopping time $\tau$,*

$$\liminf_{x\to\infty} \frac{\overline{F^{*\tau}}(x)}{\overline{F}(x)} \geq \mathbf{E}\tau.$$

It follows from Lemma 3 that it is sufficient to prove the upper bound in Theorem 1. We first discuss briefly the case where the function $h$ defined in Lemma 1 may be chosen to be additionally increasing. Here the proof of the required upper bound may proceed analogously to that of Lemma 2, working with tails rather than densities. The right-hand side of (6) is replaced by

$$\frac{\int_0^\infty x^{[t]} e^{h(x)} F^{*\tau}(dx)}{\int_0^\infty x^{[t]} e^{h(x)} F(dx)}$$

which, after integration by parts, is equal to

$$\frac{\int_0^\infty \overline{F^{*\tau}}(x)\, d(x^{[t]} e^{h(x)})}{\int_0^\infty \overline{F}(x)\, d(x^{[t]} e^{h(x)})}.$$

Thus, as in the proof of Lemma 2, we make the contrary assumption that there exist $\delta > 0$ and $x_0$ such that

$$\overline{F^{*\tau}}(x) \geq (\mathbf{E}\tau + \delta)\overline{F}(x) \qquad \text{for all } x > x_0.$$

The argument leading to a contradiction now proceeds as in the earlier proof: for the increasing function $h(x)$, as $t \to \infty$,

$$\int_0^\infty \overline{F^{*\tau}}(x)\, d(x^{[t]} e^{h(x)}) \geq (\mathbf{E}\tau + \delta + o(1)) \int_0^\infty \overline{F}(x)\, d(x^{[t]} e^{h(x)})$$

$$= (\mathbf{E}\tau + \delta + o(1)) \int_0^\infty x^{[t]} e^{h(x)} F(dx).$$

However, it is not clear that the function introduced in the statement of Lemma 1 may always be chosen to be monotone—the function constructed in the proof does not possess this property. Therefore we now use a different and novel approach which starts from the observation that the integrated tail distribution $F_I$ (see below) has a density with respect to Lebesgue measure. We apply Lemma 2 to that density. Then, in Lemma 5, we show how to use the properties of the density of $F_I$ in order to prove Theorem 1. For that to work, we restrict our consideration to the case $\tau \equiv 2$.

DEFINITION 1. For any distribution $F$ on $[0, \infty)$ with finite mean $a$, define the *integrated tail distribution* $F_I$ by

$$F_I(B) = \frac{1}{a} \int_B \overline{F}(x)\, dx.$$

The distribution $F_I$ has the decreasing density $\overline{F}(x)/a$. It is clear that both the distributions $F$ and $F_I$ are either heavy-tailed or not together. The next lemma is about the lower limit for the convolution of densities of integrated tail distributions.

LEMMA 4. *For any heavy-tailed distribution $F$ with $a \in (0, \infty]$,*

$$\liminf_{x \to \infty} \frac{1}{\overline{F}(x)} \int_0^x \overline{F}(x - y) \overline{F}(y)\, dy = 2a.$$

PROOF. First, the "lim inf" is not smaller than $2a$, because the monotonicity of the distribution tail implies the inequality

$$\frac{1}{\overline{F}(x)} \int_0^x \overline{F}(x - y) \overline{F}(y)\, dy \geq 2 \int_0^{x/2} \overline{F}(y)\, dy \to 2a.$$

Second, in the case $a < \infty$, applying Lemma 2 to the integrated tail distribution $F_I$ whose density with respect to Lebesgue measure is equal to $\overline{F}(x)/a$, we get in the special case $\tau \equiv 2$ that the "lim inf" is at most $a\mathbf{E}\tau = 2a$. The proof is complete. □

LEMMA 5. *For any heavy-tailed distribution $F$,*

$$\liminf_{x \to \infty} \frac{\overline{F * F}(x)}{\overline{F}(x)} \leq 2.$$

PROOF. We start with the case of infinite mean, that is $a = \infty$. Let $\xi_1$ and $\xi_2$ be two independent random variables with common distribution $F$. For any positive $t$, since $(\xi_1 + \xi_2)^{[t]} \leq \xi_1^{[t]} + \xi_2^{[t]}$,

$$(7) \qquad 2 = \frac{\mathbf{E}\xi_1^{[t]} + \mathbf{E}\xi_2^{[t]}}{\mathbf{E}\xi_1^{[t]}} \geq \frac{\mathbf{E}(\xi_1 + \xi_2)^{[t]}}{\mathbf{E}\xi_1^{[t]}} = \frac{\int_0^t \overline{F * F}(y)\, dy}{\int_0^t \overline{F}(y)\, dy}.$$

Suppose, contrary to the assertion of the lemma, that

$$(8) \qquad \liminf_{x \to \infty} \frac{\overline{F * F}(x)}{\overline{F}(x)} > 2.$$

Since $a = \infty$, $\int_0^t \overline{F}(y)\, dy \to \infty$ as $t \to \infty$. It then follows from the assumption (8) that

$$\liminf_{t \to \infty} \frac{\int_0^t \overline{F * F}(y)\, dy}{\int_0^t \overline{F}(y)\, dy} > 2,$$

which contradicts (7). This completes the proof in the case $a = \infty$ (a more complicated variant of the proof in this case may be found in [12]).

Now suppose $a < \infty$. Consider the distribution $G$ defined by its tail $\overline{G}(x) = (\overline{F}(x-1) + \overline{F}(x))/2$. Let $b$ denote the mean of $G$. By Lemma 4, for some $x_n \to \infty$,

$$(9) \qquad \lim_{n \to \infty} \frac{1}{\overline{G}(x_n)} \int_0^{x_n} \overline{G}(x_n - y)\overline{G}(y)\, dy = 2b.$$

For any fixed positive $t$,

$$\int_0^x \overline{G}(x-y)\overline{G}(y)\, dy$$

$$= 2 \int_0^{x/2} \overline{G}(x-y)\overline{G}(y)\, dy$$

$$\geq 2\overline{G}(x) \int_0^{x/2} \overline{G}(y)\, dy + 2(\overline{G}(x-t) - \overline{G}(x)) \int_t^{x/2} \overline{G}(y)\, dy$$

$$\sim 2\overline{G}(x)b + 2(\overline{G}(x-t) - \overline{G}(x)) \int_t^{\infty} \overline{G}(y)\, dy$$

as $x \to \infty$. Then, by (9), $\overline{G}(x_n - t) \sim \overline{G}(x_n)$ as $n \to \infty$. Equivalently, for any fixed integer $t \geq 1$,

$$\overline{F}(x_n - t - 1) + \overline{F}(x_n - t) \sim \overline{F}(x_n - t) + \overline{F}(x_n - t + 1).$$

This implies the equivalence $\overline{F}(x_n - t - 1) \sim \overline{F}(x_n - t + 1)$ and, therefore,

$$\overline{F}(x_n - t - 1) \sim \overline{F}(x_n).$$

It follows that we can choose a sufficiently slowly increasing integer sequence $t_n \to \infty$ such that $2t_n \leq x_n$ and

$$(10) \qquad \overline{F}(x_n - t_n - 1) \sim \overline{F}(x_n) \qquad \text{as } n \to \infty.$$

Indeed, there exists an increasing sequence $\{N(k)\}$ such that, for any $k \in \mathbf{N}$,

$$\frac{\overline{F}(x_n - k - 1)}{\overline{F}(x_n)} \leq 1 + \frac{1}{k}$$

for all $n \geq N(k)$. Now let $t_n = \min(k, [x_n/2])$ for $N(k) \leq n < N(k+1)$.

It follows from (10) that, as $n \to \infty$,

$$\int_0^{t_n} \overline{G}(x_n - y)\overline{G}(y)\, dy = \int_{x_n - t_n}^{x_n} \overline{G}(x_n - y)\overline{G}(y)\, dy \sim \overline{G}(x_n)b.$$

Together with (9) this implies the relation

$$\int_{t_n}^{x_n - t_n} \overline{G}(x_n - y)\overline{G}(y)\, dy = o(\overline{G}(x_n)) = o(\overline{F}(x_n)).$$

In particular,

$$(11) \qquad \int_{t_n}^{x_n - t_n} \overline{F}(x_n - y - 1)\overline{F}(y - 1)\, dy = o(\overline{F}(x_n))$$

as $n \to \infty$, due to (10). By (10) we have as well

$$(12) \qquad \int_{x_n - t_n}^{x_n} \overline{F}(x_n - y)F(dy) \leq F(x_n - t_n, x_n] = o(\overline{F}(x_n)).$$

The inner integral in the convolution formula for $F * F$ can be estimated in the following way:

$$
\begin{aligned}
\int_{t_n}^{x_n - t_n} \overline{F}(x_n - y)F(dy) &\leq \sum_{k=t_n}^{[x_n - t_n]} \int_k^{k+1} \overline{F}(x_n - y)F(dy) \\
&\leq \sum_{k=t_n}^{[x_n - t_n]} \overline{F}(x_n - k - 1)\overline{F}(k) \\
(13) \qquad &\leq \sum_{k=t_n}^{[x_n - t_n]} \int_k^{k+1} \overline{F}(x_n - y - 1)\overline{F}(y - 1)\, dy \\
&\leq \int_{t_n}^{x_n - t_n + 1} \overline{F}(x_n - y - 1)\overline{F}(y - 1)\, dy \\
&= o(\overline{F}(x_n)),
\end{aligned}
$$

by the estimate (11). It follows from (12), (13) and (10) that

$$\int_0^{x_n} \overline{F}(x_n - y) F(dy) = \int_0^{t_n} \overline{F}(x_n - y) F(dy) + o(\overline{F}(x_n))$$

$$\sim \overline{F}(x_n) \qquad \text{as } n \to \infty.$$

Hence, $\overline{F * F}(x_n) \sim 2\overline{F}(x_n)$ which concludes the proof of Lemma 5.  $\square$

**5. Light tails: local properties.** Let $\mu$ be either Lebesgue measure on $[0, \infty)$ or the counting measure on nonnegative integers. Suppose that the distribution $F$ on $[0, \infty)$ is absolutely continuous with respect to $\mu$, and let $f(x)$ be the corresponding density. Let $\tau$ be an independent stopping time. Then the distribution of the sum $S_\tau = \xi_1 + \cdots + \xi_\tau$ has density $f^{*\tau}(x)$ with respect to $\mu$.

LEMMA 6.   *If $0 < \hat{\gamma} < \infty$ and $\varphi(\hat{\gamma}) < \infty$, then*

$$\liminf_{x \to \infty} \frac{f^{*\tau}(x)}{f(x)} \leq \mathbf{E}\tau \varphi^{\tau-1}(\hat{\gamma}),$$

*provided $\mathbf{E}e^{(\hat{\gamma}+\kappa)\tau} < \infty$ for some $\kappa > 0$.*

PROOF.   Apply the exponential change of measure with parameter $\hat{\gamma}$ and consider the distribution $G$ with the density $g(x) = e^{\hat{\gamma}x} f(x)/\varphi(\hat{\gamma})$. Consider also an independent random variable $\eta$ with the distribution

$$\mathbf{P}\{\eta = n\} = \varphi^n(\hat{\gamma})\mathbf{P}\{\tau = n\}/\mathbf{E}\varphi^\tau(\hat{\gamma}).$$

The density $g^{*\eta}(x)$ of the distribution $G^{*\eta}$ is equal to

$$\sum_{n=1}^\infty g^{*n}(x)\mathbf{P}\{\eta = n\} = \frac{1}{\mathbf{E}\varphi^\tau(\hat{\gamma})} \sum_{n=1}^\infty e^{\hat{\gamma}x} f^{*n}(x)\mathbf{P}\{\tau = n\} = \frac{e^{\hat{\gamma}x} f^{*\tau}(x)}{\mathbf{E}\varphi^\tau(\hat{\gamma})}.$$

It follows from the definition of $\hat{\gamma}$ that the distribution $G$ is heavy-tailed. In addition, $\mathbf{E}e^{\kappa\eta} < \infty$. Hence, by Lemma 2,

$$\liminf_{x \to \infty} \frac{g^{*\eta}(x)}{g(x)} \leq \mathbf{E}\eta = \frac{\mathbf{E}\tau \varphi^\tau(\hat{\gamma})}{\mathbf{E}\varphi^\tau(\hat{\gamma})}.$$

Therefore,

$$\liminf_{x \to \infty} \frac{f^{*\tau}(x)}{f(x)} = \mathbf{E}\varphi^{\tau-1}(\hat{\gamma}) \liminf_{x \to \infty} \frac{g^{*\eta}(x)}{g(x)} \leq \mathbf{E}\tau \varphi^{\tau-1}(\hat{\gamma}),$$

which completes the proof.  $\square$

**6. Light tails: proof of Theorem 2.** We start with the lower bound.

LEMMA 7. *Let $\gamma$ be a positive number. If, for any fixed $y$,*

$$\liminf_{x \to \infty} \frac{\overline{F}(x-y)}{\overline{F}(x)} \geq e^{\gamma y}, \tag{14}$$

*then*

$$\liminf_{x \to \infty} \frac{\overline{F * F}(x)}{\overline{F}(x)} \geq 2\varphi(\gamma);$$

*the case $\varphi(\gamma) = \infty$ is not excluded.*

PROOF. For $x > 2t$, we have the inequality

$$\overline{F * F}(x) \geq 2 \int_0^t \overline{F}(x-y)F(dy) = 2\overline{F}(x) \int_0^t \frac{\overline{F}(x-y)}{\overline{F}(x)} F(dy). \tag{15}$$

Now by condition (14) and Fatou's lemma,

$$\liminf_{x \to \infty} \frac{\overline{F * F}(x)}{\overline{F}(x)} \geq 2 \int_0^t \liminf_{x \to \infty} \frac{\overline{F}(x-y)}{\overline{F}(x)} F(dy) \geq 2 \int_0^t e^{\gamma y} F(dy).$$

Letting $t \to \infty$ we arrive at conclusion of the lemma. □

The following auxiliary lemma compares the tail behavior of the convolution tail and that of the exponentially transformed distribution.

LEMMA 8. *Let the distribution $F$ and the number $\gamma \geq 0$ be such that $\varphi(\gamma) < \infty$. Let the distribution $G$ be the result of the exponential change of measure with parameter $\gamma$, that is, $G(du) = e^{\gamma u}F(du)/\varphi(\gamma)$. Then*

$$\liminf_{x \to \infty} \frac{\overline{G * G}(x)}{\overline{G}(x)} \geq \frac{1}{\varphi(\gamma)} \liminf_{x \to \infty} \frac{\overline{F * F}(x)}{\overline{F}(x)}$$

*and*

$$\limsup_{x \to \infty} \frac{\overline{G * G}(x)}{\overline{G}(x)} \leq \frac{1}{\varphi(\gamma)} \limsup_{x \to \infty} \frac{\overline{F * F}(x)}{\overline{F}(x)}.$$

PROOF. Put

$$\hat{c} \equiv \liminf_{x \to \infty} \frac{\overline{F * F}(x)}{\overline{F}(x)}.$$

By Lemma 3, $\hat{c} \in [2, \infty]$. For any fixed $c \in (0, \hat{c})$, there exists $x_0 > 0$ such that, for any $x > x_0$,

$$\overline{F * F}(x) \geq c\overline{F}(x). \tag{16}$$

Integration by parts yields

$$(17) \quad \overline{G * G}(x) = \frac{1}{\varphi^2(\gamma)} \int_x^\infty e^{\gamma y}(F * F)(dy)$$
$$= \frac{1}{\varphi^2(\gamma)} \left[ e^{\gamma x} \overline{F * F}(x) + \int_x^\infty \overline{F * F}(y) \, de^{\gamma y} \right].$$

Using also (16) we get, for $x > x_0$,

$$\overline{G * G}(x) \geq \frac{c}{\varphi^2(\gamma)} \left[ e^{\gamma x} \overline{F}(x) + \int_x^\infty \overline{F}(y) \, de^{\gamma y} \right]$$
$$= \frac{c}{\varphi^2(\gamma)} \int_x^\infty e^{\gamma y} F(dy) = \frac{c}{\varphi(\gamma)} \overline{G}(x).$$

Letting $c \uparrow \hat{c}$, we obtain the first conclusion of the lemma. The proof of the second conclusion follows similarly. $\quad\square$

REMARK 1. For $\gamma < 0$, the statements of Lemma 8 need not hold. The reason is that, for $\gamma < 0$, the function $e^{\gamma y}$ is decreasing and while the first term in brackets in (17) is positive the second one is negative.

LEMMA 9. If $0 < \hat{\gamma} < \infty$ and $\varphi(\hat{\gamma}) < \infty$, then

$$\liminf_{x \to \infty} \frac{\overline{F * F}(x)}{\overline{F}(x)} \leq 2\varphi(\hat{\gamma})$$

and

$$\limsup_{x \to \infty} \frac{\overline{F * F}(x)}{\overline{F}(x)} \geq 2\varphi(\hat{\gamma}).$$

PROOF. We apply the exponential change of measure with parameter $\hat{\gamma}$ and consider the distribution $G(du) = e^{\hat{\gamma} u} F(du)/\varphi(\hat{\gamma})$. Again from the definition of $\hat{\gamma}$, the distribution $G$ is heavy-tailed. Hence,

$$\limsup_{x \to \infty} \frac{\overline{G * G}(x)}{\overline{G}(x)} \geq \liminf_{x \to \infty} \frac{\overline{G * G}(x)}{\overline{G}(x)} = 2,$$

by Theorem 1. The result now follows from Lemma 8 with $\gamma = \hat{\gamma}$. $\quad\square$

**7. Convolution tail equivalent distributions: proof of Theorem 3.** In the case where $F$ is heavy-tailed we have $\hat{\gamma} = 0$ and $\varphi(\hat{\gamma}) = 1$. By Theorem 1, $c = 2$ as required.

In the case $\hat{\gamma} \in (0, \infty)$ and $\varphi(\hat{\gamma}) < \infty$, the desired conclusion follows from Lemma 9.

Finally, consider the case $\varphi(\hat\gamma) = \infty$. Fix an arbitrary $\gamma \in (0, \hat\gamma)$. Then $\varphi(\gamma) < \infty$. Consider the distribution $G(du) = e^{\gamma u} F(du)/\varphi(\gamma)$. Then Lemma 8 shows that

$$\overline{G * G}(x) \sim \frac{c}{\varphi(\gamma)} \overline{G}(y).$$

This equivalence and Lemma 3 imply the inequality $c/\varphi(\gamma) \geq 2$. Since $\varphi(\gamma) \uparrow \varphi(\hat\gamma) = \infty$ as $\gamma \uparrow \hat\gamma$, we obtain $c = \infty$. The proof is complete.

**8. Corollaries for the classes $\mathcal{S}$, $\mathcal{S}(\gamma)$ and $\mathcal{S}^*$.** In this section we continue to consider distributions on $[0, \infty)$ only.

DEFINITION 2. A distribution $F$ with unbounded support is called *long-tailed* ($F \in \mathcal{L}$) if, for any fixed $y$, $\overline{F}(x + y) \sim \overline{F}(x)$ as $x \to \infty$. Clearly, any long-tailed distribution is heavy-tailed.

DEFINITION 3. A distribution $F$ with unbounded support belongs to the class $\mathcal{S}$ of *subexponential distributions* if $\overline{F * F}(x) \sim 2\overline{F}(x)$ as $x \to \infty$. It is known from [1] that $\mathcal{S} \subset \mathcal{L}$.

DEFINITION 4. A distribution $F$ with finite mean $a$ belongs to the class $\mathcal{S}^*$ if

$$\int_0^x \overline{F}(x - y)\overline{F}(y) \, dy \sim 2a\overline{F}(x) \qquad \text{as } x \to \infty.$$

It is known (see [8]) that $F \in \mathcal{S}^*$ implies $F \in \mathcal{S}$ and $F^I \in \mathcal{S}$. The converse implication is not true, in general; see [5].

Theorem 1 implies the following result related to the definition of $\mathcal{S}$.

THEOREM 4. *Let the distribution $F$ on $[0, \infty)$ be heavy-tailed. If, for some $c \in (0, \infty)$, $\overline{F * F}(x) \sim c\overline{F}(x)$ as $x \to \infty$, then $F \in \mathcal{S}$.*

Theorem 4 generalizes the main theorem of [10] where it was additionally assumed that $F$ is long-tailed. This result, for long-tailed distribution, was first formulated by Chover, Ney and Wainger ([3], page 664); the corresponding proof, based on Banach algebra techniques, contains some holes; see the comments by Rogozin and Sgibnev ([11], Section 4) on the matter. A further attempted proof by Cline ([4], Theorem 2.9) also contains a gap [in the proof of Lemma 2.3(ii); in particular, it was not proved in line $-7$ on page 351 that one can choose $t_0$ independently of $n$]. To the best of our knowledge, the only paper which states the same result as Theorem 4 is that of Teugels ([13], Theorem 1(i)), but the proof there is incorrect in line $-5$ on page 1002 and in lines 11–12 on page 1003.

From Lemma 4 we get the following result for the class $\mathcal{S}^*$.

THEOREM 5. *Let the distribution $F$ on $[0,\infty)$ be heavy-tailed. If, for some $c \in (0,\infty)$,*

$$\int_0^x \overline{F}(x-y)\overline{F}(y)\,dy \sim c\overline{F}(x) \qquad as\ x \to \infty,$$

*then $F$ has a finite mean and $F \in \mathcal{S}^*$.*

DEFINITION 5. A distribution $F$ with unbounded support belongs to the class $\mathcal{S}(\gamma)$, $\gamma \geq 0$, if the following conditions hold:

  (i) $\varphi(\gamma) < \infty$;
 (ii) for any fixed $y$, $\overline{F}(x+y)/\overline{F}(x) \to e^{-\gamma y}$ as $x \to \infty$;
(iii) $\overline{F*F}(x) \sim 2\varphi(\gamma)\overline{F}(x)$ as $x \to \infty$.

It follows that if $F \in \mathcal{S}(\gamma)$, then $\hat{\gamma} = \gamma < \infty$. Further, the class $\mathcal{S}(0)$ coincides with the class $\mathcal{S}$; see Definition 3 above.

In the following theorem we observe a lifting property for the class $\mathcal{S}(\gamma)$, which for the case $\beta = \gamma$ was originally proved by Embrechts and Goldie ([7], Theorem 3.1).

THEOREM 6. *If $F \in \mathcal{S}(\hat{\gamma})$, then, for any $\beta \in [0,\hat{\gamma}]$, the distribution $G(du) = e^{\beta u}F(du)/\varphi(\beta) \in \mathcal{S}(\hat{\gamma}-\beta)$. In particular, the distribution $G(du) = e^{\hat{\gamma}u}F(du)/\varphi(\hat{\gamma})$ is subexponential.*

PROOF. First note that the inverse implication does not in general hold (see [7], Theorem 3.1). The main reason for this is pointed out in Remark 1 above.

The Laplace transform of $G$ at the point $\hat{\gamma} - \beta$ is equal to $\varphi(\hat{\gamma})/\varphi(\beta)$. By Lemma 8, the distribution $G$ satisfies property (iii) of Definition 5 with $\gamma$ replaced by $\hat{\gamma} - \beta$. For the case $\beta = \hat{\gamma}$ this completes the proof. Now consider the case $\beta < \hat{\gamma}$: we require to prove (ii). Integration by parts yields

$$\begin{aligned}
e^{(\hat{\gamma}-\beta)x}&\overline{G}(x)\\
&= \frac{1}{\varphi(\beta)}\left[e^{\hat{\gamma}x}\overline{F}(x) + e^{(\hat{\gamma}-\beta)x}\beta\int_x^\infty \overline{F}(y)e^{\beta y}\,dy\right]\\
&= \frac{1}{\varphi(\beta)}\left[e^{\hat{\gamma}x}\overline{F}(x) + e^{(\hat{\gamma}-\beta)x}\beta\overline{F}(x)e^{\beta x}\int_x^\infty \frac{\overline{F}(y)e^{\hat{\gamma}y}}{\overline{F}(x)e^{\hat{\gamma}x}}e^{(\beta-\hat{\gamma})(y-x)}\,dy\right]\\
&\sim \frac{1}{\varphi(\beta)}\left[e^{\hat{\gamma}x}\overline{F}(x) + e^{\hat{\gamma}x}\beta\overline{F}(x)\int_x^\infty e^{(\beta-\hat{\gamma})(y-x)}\,dy\right]
\end{aligned}$$

as $x \to \infty$, since $\beta - \hat{\gamma} < 0$ and by the dominated convergence theorem. Hence,

$$e^{(\hat{\gamma}-\beta)x}\overline{G}(x) \sim \frac{1}{\varphi(\beta)}\left[e^{\hat{\gamma}x}\overline{F}(x) + e^{\hat{\gamma}x}\frac{\beta}{\hat{\gamma}-\beta}\overline{F}(x)\right] = \frac{1}{\varphi(\beta)}e^{\hat{\gamma}x}\overline{F}(x)\frac{\hat{\gamma}}{\hat{\gamma}-\beta},$$

which implies that for any fixed $t$, $e^{(\hat{\gamma}-\beta)x}\overline{G}(x) \sim e^{(\hat{\gamma}-\beta)(x+t)}\overline{G}(x+t)$ as $x \to \infty$. The proof is complete.  □

We would like to formulate the following hypothetically equivalent definition of the class $\mathcal{S}(\gamma)$ which instead of properties (i) and (ii) from Definition 5 assumes only that $\gamma$ is the right convergence point of the Laplace transform.

CONJECTURE 1.  *The distribution $F$ on $[0, \infty)$ with unbounded support belongs to the class $\mathcal{S}(\gamma)$ if and only if:*

(a)  $\gamma = \hat{\gamma}$;
(b)  *for some $c \in [2, \infty)$, $\overline{F * F}(x) \sim c\overline{F}(x)$ as $x \to \infty$.*

For $\gamma > 0$ we have neither a proof of the conjecture nor a counterexample. We can prove only the following weakened version of this statement.

THEOREM 7.  *The distribution $F$ on $[0, \infty)$ with unbounded support belongs to the class $\mathcal{S}(\gamma)$ if and only if:*

(a)  *and* (b) *of Conjecture 1 hold;*
(c)  *condition* (14) *holds.*

Theorem 7, with the condition (c) replaced by the stronger requirement

$$\frac{\overline{F}(x-y)}{\overline{F}(x)} \to e^{\gamma x} \qquad \text{for all } y > 0,$$

was proved by Rogozin and Sgibnev [11]. The proofs of earlier assertions of this latter result by Chover, Ney and Wainger [3] and by Cline [4] are incomplete for the reasons mentioned after Theorem 4 above. Cline's version of the result [4] was referenced by Pakes [9] in his study of distributions on the whole real line. Theorem 7 is close in spirit to an assertion of Teugels ([13], Theorem 1(ii)). However, Embrechts and Goldie ([7], Section 3) showed that this assertion is incorrect.

PROOF OF THEOREM 7.  By Theorem 3, $c = 2\varphi(\gamma) < \infty$. Thus it suffices to prove that the condition (ii) of Definition 5 is satisfied. Suppose that, on the contrary, there exist $y > 0$, $\delta > 0$, and a sequence $x_n \to \infty$ such that, for any $n \geq 1$,

$$\frac{\overline{F}(x_n - y)}{\overline{F}(x_n)} \geq e^{\gamma y} + 3\delta.$$

Since the function $\overline{F}(x_n - y)$ is increasing in $y$, there exists $y_0 > y$ such that, for any $n \geq 1$ and $t \in [y, y_0]$,

$$\frac{\overline{F}(x_n - t)}{\overline{F}(x_n)} \geq e^{\gamma t} + 2\delta.$$

Choose any $z$ so that $F(z + y, z + y_0] \equiv \overline{F}(z + y) - \overline{F}(z + y_0) > 0$. Since

$$\liminf_{n \to \infty} \frac{\overline{F}(x_n)}{\overline{F}(x_n + z)} \geq e^{\gamma z}$$

by the condition (14), for all sufficiently large $n$ and $t \in [z + y, z + y_0]$,

$$\frac{\overline{F}(x_n + z - t)}{\overline{F}(x_n + z)} = \frac{\overline{F}(x_n + z - t)}{\overline{F}(x_n)} \frac{\overline{F}(x_n)}{\overline{F}(x_n + z)} \geq e^{\gamma t} + \delta.$$

Together with (14), (15) and Fatou's lemma this implies

$$\liminf_{n \to \infty} \frac{\overline{F * F}(x_n + z)}{\overline{F}(x_n + z)} \geq 2\varphi(\gamma) + 2\delta F[z + y, z + y_0].$$

The latter sum exceeds $2\varphi(\gamma)$. This contradicts the equality $c = 2\varphi(\gamma)$ and thus completes the proof. $\square$

We now consider the local properties of convolutions. For simplicity we study the lattice case only. Hereinafter denote by $F\{n\}$ the distribution mass at the point $n$.

DEFINITION 6. A distribution $F$ on $\mathbf{Z}^+$ with unbounded support belongs to the class $\mathcal{S}_{\text{lattice}}(\gamma)$, $\gamma \geq 0$, if the following conditions hold:

(i)  $\varphi(\gamma) < \infty$;
(ii)  $F\{n + 1\}/F\{n\} \to e^{-\gamma}$ as $n \to \infty$;
(iii)  $F * F\{n\} \sim 2\varphi(\gamma)F\{n\}$ as $n \to \infty$.

THEOREM 8. A lattice distribution $F$ belongs to the class $\mathcal{S}_{\text{lattice}}(\hat{\gamma})$ if and only if:

(i)  $\liminf_{x \to \infty} F\{n - 1\}/F\{n\} \geq e^{\hat{\gamma}}$;
(ii)  $F * F\{n\} \sim cF\{n\}$ as $n \to \infty$ for some $c \in (0, \infty)$.

PROOF. This theorem generalizes the corresponding assertions by Chover, Ney and Wainger [2] and by Embrechts ([6], Theorem 2.8). To prove Theorem 8 first note that $c = 2\varphi(\hat{\gamma})$ with necessity and this fact correlates with Theorem 3. Indeed, Lemmas 2 and 6 with $\tau \equiv 2$ and counting measure for $\mu$ imply $c \leq 2\varphi(\hat{\gamma})$. Further, the estimate

$$\frac{F * F\{n\}}{F\{n\}} = \sum_{k=0}^{n} \frac{F\{k\}F\{n - k\}}{F\{n\}} \geq 2 \sum_{k=0}^{[n/2]} F\{k\} \frac{F\{n - k\}}{F\{n\}}$$

implies $c \geq 2\varphi(\hat{\gamma})$, due to condition (i) of the theorem. Now it suffices to establish the condition (ii) of Definition 6. Suppose that, on the contrary, there exist $\delta > 0$ and a sequence $n_i \to \infty$ such that, for any $i \geq 1$,

$$\frac{F\{n_i - 1\}}{F\{n_i\}} \geq e^{\gamma y} + \delta.$$

Then arguments similar to those of the proof of Theorem 7 lead to $c > 2\varphi(\gamma)$. This contradiction completes the proof.   $\square$

**9. Counterexamples.**   Let $F$ be an atomic distribution at the points $x_n$, $n = 0, 1, \ldots$, with masses $p_n$, that is, $F\{x_n\} = p_n$. Suppose that $x_0 = 1$ and that $x_{n+1} > 2x_n$ for every $n$. Then the tail of the convolution $F * F$ at the point $x_n - 1$ is equal to

$$\overline{F*F}(x_n - 1) = (F \times F)([x_n, \infty) \times [0, \infty)) + (F \times F)([0, x_{n-1}) \times [x_n, \infty))$$

$$\sim 2\overline{F}(x_n - 1) \qquad \text{as } n \to \infty.$$

Hence,

$$\lim_{n \to \infty} \frac{\overline{F*F}(x_n - 1)}{\overline{F}(x_n - 1)} = 2.$$

From this equality and Lemma 3

$$(18) \qquad\qquad \liminf_{x \to \infty} \frac{\overline{F*F}(x)}{\overline{F}(x)} = 2.$$

We now consider three particular cases.

EXAMPLE 1.   Let $x_n = 3^n$, $n = 0, 1, \ldots$, and $p_n = ce^{-\hat{\gamma} 3^n} 3^{-n}$, where $\hat{\gamma} > 0$ and $c$ is the normalizing constant. We have $\varphi(\hat{\gamma}) = 3c/2 \in (1, \infty)$ and $\varphi(\gamma) = \infty$ for any $\gamma > \hat{\gamma}$. In this example

$$\liminf_{x \to \infty} \frac{\overline{F*F}(x)}{\overline{F}(x)} = 2 < 2\varphi(\hat{\gamma}).$$

Clearly, the condition (2) is violated for this distribution.

EXAMPLE 2.   Let $x_n = 3^n$, $n = 0, 1, \ldots$, and $p_n = ce^{-\hat{\gamma} 3^n}$. In this example we have $\varphi(\gamma) < \infty$ for any $\gamma < \hat{\gamma}$ and $\varphi(\hat{\gamma}) = \infty$. Nevertheless, the relation (18) still holds.

EXAMPLE 3.   Let $p_n = ce^{-x_n^2}$. In this example we have $\varphi(\gamma) < \infty$ for any $\gamma$, that is, $\hat{\gamma} = \infty$. But again the relation (18) remains valid.

**10. Convolutions of nonidentical distributions.**

THEOREM 9. *Let $F_1$ and $F_2$ be two distributions on $[0, \infty)$. If the distribution $F_1$ is heavy-tailed, then*

$$\liminf_{x \to \infty} \frac{\overline{F_1 * F_2}(x)}{\overline{F}_1(x) + \overline{F}_2(x)} = 1. \tag{19}$$

PROOF. For any two distributions $F_1$ and $F_2$ on $[0, \infty)$,

$$\overline{F_1 * F_2}(x) \geq (F_1 \times F_2)((x, \infty) \times [0, x]) + (F_1 \times F_2)([0, x] \times (x, \infty))$$
$$= \overline{F}_1(x)F_2(x) + F_1(x)\overline{F}_2(x) \tag{20}$$
$$\sim \overline{F}_1(x) + \overline{F}_2(x) \qquad \text{as } x \to \infty,$$

which implies that the left-hand side of (19) is at least 1. Assume now that it is strictly greater than 1. Then there exists $\varepsilon > 0$ such that, for all sufficiently large $x$,

$$\frac{\overline{F_1 * F_2}(x)}{\overline{F}_1(x) + \overline{F}_2(x)} \geq 1 + 2\varepsilon. \tag{21}$$

Consider the distribution $G = (F_1 + F_2)/2$. This distribution is heavy-tailed. By Theorem 1 we get

$$\liminf_{x \to \infty} \frac{\overline{G * G}(x)}{\overline{G}(x)} = 2. \tag{22}$$

On the other hand, (21) and (20) imply that, for all sufficiently large $x$,

$$\overline{G * G}(x) = \frac{\overline{F_1 * F_1}(x) + \overline{F_2 * F_2}(x) + 2\overline{F_1 * F_2}(x)}{4}$$
$$\geq \frac{2(1 - \varepsilon)\overline{F}_1(x) + 2(1 - \varepsilon)\overline{F}_2(x) + 2(1 + 2\varepsilon)(\overline{F}_1(x) + \overline{F}_2(x))}{4}$$
$$= 2(1 + \varepsilon/2)\overline{G}(x),$$

which contradicts (22). The proof is complete. $\square$

**Acknowledgments.** The final version of this paper was written while both authors were visiting the Boole Centre for Research in Informatics, University College Cork, thanks to the hospitality of Neil O'Connell and financial support of Science Foundation Ireland. We are grateful to Stan Zachary for extensive help which led to a significant improvement of the style and clearness of exposition. We would like to thank the anonymous referee for his/her very positive comments and suggestions.

HERIOT-WATT UNIVERSITY
EDINBURGH EH14 4AS
UNITED KINGDOM
E-MAIL: S.Foss@ma.hw.ac.uk

SOBOLEV INSTITUTE OF MATHEMATICS
4 KOPTYUGA PR.
NOVOSIBIRSK 630090
RUSSIA
E-MAIL: korshunov@math.nsc.ru